\documentclass[reqno,12pt,english]{amsart}
\usepackage{amssymb,amsmath,amstext,amsthm,amsfonts,amscd}
\usepackage[dvips]{graphicx}
\usepackage{color}
\usepackage{psfrag,euscript,a4wide,stmaryrd,wasysym}

\newtheorem{theorem}{Theorem}[section]
\newtheorem{corollary}{Corollary}[section]

\newtheorem{maintheorem}{Theorem}

\newcommand{\cmt}{\begin{maintheorem}}
\newcommand{\fmt}{\end{maintheorem}}

\newtheorem{maincorollary}[maintheorem]{Corollary}

\newcommand{\cmc}{\begin{maincorollary}}
\newcommand{\fmc}{\end{maincorollary}}

\newtheorem{T}{Theorem}[section]
\newcommand{\cte}{\begin{T}}
\newcommand{\fte}{\end{T}}

\newtheorem{Corollary}[T]{Corollary}
\newcommand{\cco}{\begin{Corollary}}
\newcommand{\fco}{\end{Corollary}}

\newtheorem{Proposition}[T]{Proposition}
\newcommand{\cpr}{\begin{Proposition}}
\newcommand{\fpr}{\end{Proposition}}

\newtheorem{Lemma}[T]{Lemma}
\newcommand{\cle}{\begin{Lemma}}
\newcommand{\fle}{\end{Lemma}}

\newcommand{\csle}{\begin{Lemma}}
\newcommand{\fsle}{\end{Lemma}}

\newtheorem{Remark}[T]{Remark}
\newcommand{\cre}{\begin{Remark}}
\newcommand{\fre}{\end{Remark}}

\newtheorem{Definition}[T]{Definition}
\newcommand{\cde}{\begin{Definition}}
\newcommand{\fde}{\end{Definition}}

\newcommand{\SA}{{\mathcal A}}

\newcommand{\SC}{{\mathcal C}}

\newcommand{\SF}{{\mathcal F}}

\newcommand{\SK}{{\mathcal K}}

\newcommand{\SM}{{\mathcal M}}

\newcommand{\SP}{{\mathcal{P}}}

\newcommand{\SR}{{\mathcal R}}
\renewcommand{\SS}{{\mathcal S}}
\newcommand{\ST}{{\mathcal{ T}}}
\newcommand{\SU}{{\mathcal U}}

\newcommand{\dem}{\begin{proof}}
\newcommand{\cqd}{\end{proof}}

\newtheorem{proposition}[theorem]{Proposition}
\newtheorem{lemma}[theorem]{Lemma}

\theoremstyle{definition}
\newtheorem{remark}[theorem]{Remark}

\newtheorem{definition}[theorem]{Definition}

\newcommand{\field}[1]{\mathbb{#1}}

\newcommand{\ov} {\overline}

\newcommand{\re}{\field{R}}

\renewcommand{\natural}{\field{N}}

    \newcommand{\Ga}{\Gamma}

\newcommand{\la} {\lambda}      \newcommand{\La}{\Lambda}

\newcommand{\vsi}{\varsigma}

\newcommand{\vr}{\varphi}

\newcommand{\NE}{\operatorname{{\text{\em NEnd}^1}}}
\newcommand{\DI}{\operatorname{{\text{\em Diff}^1}}}

\newcommand{\diam}{\operatorname{{diam}}}

\begin{document}

\author{Armando Castro}
\address{Departamento de Matem\'atica, Universidade Federal da Bahia\\
Av. Ademar de Barros s/n, 40170-110 Salvador, Brazil.}
\email{armandomat@pesquisador.cnpq.br}


\date{\today}


\thanks{Work carried out at the  Federal University of
Bahia. Partially supported by CNPq (PQ 10/2007) and UFBA}

\title[Criteria of Hyperbolicity]{Ergodic Closing and New Criteria for Hyperbolicity based on periodic Sets}

\maketitle

\begin{abstract}
We prove some criteria for uniform hyperbolicity based on the
periodic points of the transformation. More precisely, if a mild
hyperbolicity condition holds for the periodic points of any
diffeomorphism in a residual subset of a $C^1$-open set $\SU$ then
there exists an open and dense subset $\SA\subset \SU$ of Axiom A
diffeomorphisms.  Moreover, we also prove a noninvertible version of Ergodic
Closing Lemma which we use to prove a counterpart of this result for
local diffeomorphisms.
\end{abstract}

\section{Introduction}



The notion of uniform hyperbolicity was coined in the mid sixties by
the pioneering works of Smale and Anosov and constitutes a rich
class of dynamical systems.  Indeed, many geometrical, topological
and ergodic properties have been proved to hold for uniformly
hyperbolic dynamical systems in both discrete and continuous time
setting.
From the eminently topological point of view, one of the greatest
goals achieved by mathematicians in describing uniformly hyperbolic
behavior was the proof of the $C^1$-Structural Stability Conjecture,
stated by Palis and Smale in \cite{PS70}. Roughly,  a diffeomorphism
$f$ is conjugated to every diffeomorphism in a small
$C^1$-neighbourhood of it if and only if $f$ is hyperbolic and
satisfies the strong transversality condition. The contributions of
De Melo \cite{DeMelo}, Franks \cite{F},  Robbin \cite{Ro71},
Robinson\cite{Ros71}\cite{Ros76}, Pliss \cite{Pli},
Ma\~n\'e\cite{Ma86} and Liao\cite{L80}, among others, were
fundamental to characterize $C^1$-structural stability by methods
that rely on periodic orbits.

It is quite natural to understand that periodic orbits are key
ingredients to characterize uniform hyperbolicity and it is natural
to ask whether the loss of hyperbolicity can be observed at this
level. Moreover, the more recent developments on the theory of
nonuniformly hyperbolic transformations can be used to give a
positive answer to this question. In fact a recent contribution in
this direction was given by  Oliveira, Pinheiro and the present
author in \cite{COP} that,  inspired by \cite{C1,C2,Pi06}, studied
the relation between asymptotic growth rates of the periodic points
of diffeomorphisms and local diffeomorphisms on  compact manifolds,
and its relations with uniform hyperbolicity and uniform expansion.
More precisely, the authors proved that if a diffeomorphism
satisfies the shadowing property (in particular if it is conjugated
to a hyperbolic one)  and the periodic set is non-uniformly
hyperbolic (NUH)  such diffeomorphism is hyperbolic and an analogous
statement for local diffeomorphisms. One of the key ideas of
\cite{COP} is that the asymptotic hyperbolicity of a NUH periodic
set spreads out to any recurrent point using the shadowing property.
In particular, if the NUH periodic set exhibits a dominated
splitting then any $f$-invariant probability measure has only
nonzero Lyapunov exponents and, using \cite{YC},  it follows that
$f$ is uniformly hyperbolic.
However, even though dominated splitting is a generic (residual)
feature among robustly transitive maximal invariant sets  (see
\cite{CPD})  the strong notion of shadowing seems to be somewhat
rare far from uniformly hyperbolic maps as proved in \cite{DA}.
Nevertheless,  \cite{COP} gave rise to some ideas to prove uniform
hyperbolicity used e.g. in  \cite{LS08}.

In the present work  we are interested in providing new criteria to
obtain uniform hyperbolicity without assuming the strong shadowing
property. Our main assumption, in both invertible and non-invertible
settings, is just that a $C^1$-open set $\SU$ exhibits  a
residual subset of transformations whose periodic sets are NUH. Note
that, in principle, the class of transformations considered here
could be very far from the uniformly hyperbolic ones.  Nevertheless,
we are able to prove that under this hypothesis there exists a
residual subset of $\SU$ whose elements present uniform
hyperbolicity.
We refer the reader to Section~\ref{s.setting} for the definitions and precise statements.
%
%

One of the novelties of our approach is that, differently from the
strong shadowing conditions typically considered when studying the
geometrical aspects of dynamical systems, we use weak shadowing
properties that hold just for points in total probability sets. This
result, known as the Ergodic Closing Lemma, was proved by Ma\~n\'e's
\cite{Ma82} in the invertible context and it was extended by us
(see Section~\ref{s.ecl}) to the context of local diffeomorphisms.
Since this kind of shadowing by periodic points holds
residually in the space of $C^1$-diffeomophisms and local diffeomorphisms our
criteria provides a generic subset of transformations exhibiting uniform
hyperbolicity.



The paper is organized as follows. In Section~\ref{s.setting} we detail the context and
give exact statements of our theorems. The proofs of the criteria to obtain uniform hyperbolicity
from the assumption on the periodic set are given along Section~\ref{s.p.criteria} in both invertible
and noninvertible context. In the last section, we prove the version of Ergodic Closing Lemma for Endomorphisms which is needed for the proof of our Criteria in the Endomorphism case.
\vspace{.2cm}

{\bf Acknowledgements: }I am grateful to professors Paulo Varandas, Marcelo Viana, Jacob Palis Jr., Fl\'avio Abdenur, Vilton Pinheiro, and to my student Luciana
Salgado for conversations on topics related to this paper. I am also
grateful to my beloved wife Maria Teresa Gilly, whose shadow shadows mine.




 \section{Setting and Statement of Main results}\label{s.setting}

Throughout, $M$ will always denote a finite dimensional compact
Riemannian manifold. Let $\NE(M)\subset C^1(M)$  will  denote the set of $C^1$
nonsingular endomorphisms (or local diffeomorphisms) in $M$
endowed with the $C^1$ topology and
let $\DI(M)$ denote the set of $C^1$ diffeomorphisms on $M$.
First we recall some necessary definitions.

\begin{definition}
Let $\La$ be a compact invariant set for a $C^1$ diffeomorphism $f$
of a manifold $M$. We say that $\La$ is a hyperbolic set if there is
a continuous splitting  $T_\La M= E^s \oplus E^u$ which is
$Df$-invariant ($Df(E^s)= E^s, Df(E^u)= E^u$) and for which there
are constants $c> 0$, $0< \vsi< 1$, such that
$$
\|Df^n|_{E^s}\| < c \cdot \vsi^n, \quad  \|Df^{-n}|_{E^u}\| < c
\cdot \vsi^n, \forall n \in \natural.
$$
\end{definition}

\vspace{.23cm}
Recall also that a diffeomorphism $f: M \to M$ is {\em Axiom A} if
the nonwandering set $\Omega(f)$ is a hyperbolic set and
$\Omega(f)= \ov{Per(f)}$.
In such case, $\Omega(f)$ admits a decomposition
$\Omega(f)= \Omega_1 \cup \dots \cup \Omega_w$ into
closed, disjoint transitive subsets. A cycle on $\Omega(f)$ is
a sequence $\Omega_{i_1}, \dots, \Omega_{i_k}$ with points
$x_1, y_1 \in \Omega_{i_1}$, $\dots x_k, y_k \in \Omega_{i_k}$
such that $W^s(x_1) \cap W^u(y_2) \neq \emptyset$,
$\dots W^s(x_k) \cap W^u(y_1) \neq \emptyset$.




\begin{definition} \label{defNUH} 
\label{defnuh} Given $f\in \DI(M)$ we say that an $f$-invariant set $X \subset M$ is a {\em
non uniformly hyperbolic set} (or simply NUH set) if
\begin{enumerate}
\item
There is an $Df-$invariant splitting $T_X M= E^{cs} \oplus E^{cu}$;

\item There exists $\la< 0$ and an
adapted Riemannian metric for which any point $p \in X$  satisfies
$$
\liminf_{n \to +\infty}\frac{1}{n}\sum_{j= 0}^{n-1} \log
\|Df(f^j(p))|_{E^{cs}(f^j(p))}\| \leq \la
$$
and
$$
\liminf_{n \to +\infty}\frac{1}{n}\sum_{j= 0}^{n- 1} \log
\|[Df(f^{j}(p))|_{E^{cu}(f^{j}(p))}]^{-1}\| \leq \la
$$
\end{enumerate}
\end{definition}

For simplicity sometimes we will say that $f$ is NUH on $X$ meaning
that $X\subset M$ is non uniformly hyperbolic for $f$.  Note that if
a set $X$ is NUH and $X$ has positive measure for some ergodic
measure  $\mu$, then $\mu$ has only nonzero Liapunov exponents.  The
properties of such measures for $C^{1+\alpha}$ diffeomorphisms are
extensively studied in Pesin Theory (see e.g. \cite{Pe76},
\cite{PeSi82}).

\begin{definition} \label{defdom}
Given $f\in \DI(M)$ and
an invariant set $X \subset M$, we say that a $Df$-invariant splitting
 $T_X M= E \oplus \hat E$ is a {\em dominated splitting} if for some $l \geq 1$
 there exists $0<\eta<1$ so that
$$
\sup_{v \in E, \|v\|= 1} \{\|Df^l(x) v\|\} \cdot (\inf_{v \in \hat
E, \|v\|= 1} \{\|Df^l(x) v\|\})^{-1} \leq \eta, \, \forall x \in X.
$$
\end{definition}




As discussed in the introduction, this paper is devoted to study the
consequences for an open set $\SU$ of $C^1$-transformations that
exhibits a residual subset $\SS$ in which each map exhibits a NUH
periodic set admitting a dominated splitting. Our first main results
are as follows.

 \cmt \label{theo2} Let $\SU \subset \DI(M)$ be an open subset
of diffeomorphisms and assume that every $f$ in some residual subset
$\SS$ of $\SU$, the set $Per(f) \subset M$ of periodic points of $f$
is non uniformly hyperbolic (NUH), and
 $T_{Per(g)}M= E^{cs} \oplus E^{cu}$ is a dominated splitting. Then,
there exists an open and dense subset of $\SU$ whose elements are Axiom A
diffeomorphisms (with no cycles). In particular, $\SU$ is contained in the closure of
the Axiom A diffeomorphisms set.
\fmt

In fact, this result is a consequence of the more general:

\cmt \label{theo3}  Let $\SU \subset \DI(M)$ be an open subset of
diffeomorphisms. Suppose that for any $f$ in some residual subset
$\SS$  of $\SU$, the set $Per(f)$ of periodic points of $f$ is non
uniformly hyperbolic (NUH), and $T_{Per(f)}M= E^{cs} \oplus E^{cu}$
is a continuous splitting extending to $T_{\ov{Per(f)}}M$. Then,
there exists an open and dense subset $\SA$ of $\SU$ whose elements are
Axiom A diffeomorphisms (with no cycles). In particular, $\SU$ is contained in the
closure of the Axiom A diffeomorphisms set.
\fmt

We point out that although Theorem \ref{theo2} is a consequence of
Theorem \ref{theo3},  the hypotheses in Theorem \ref{theo2} are
easier to verify. An important remark is also that in principle,
due to Palis' work (see \cite{Pa70}), Axiom A are not necessarily open
if  a no-cycles condition is not assumed. We will need a periodic set semicontinuity argument
in order to obtain the no-cycles condition and the openness stated
in both theorems above.
It is also interesting to notice that, at least generically,
bifurcations from uniformly hyperbolic diffeomorphisms are verified
by some lack of hyperbolicity at periodic points. See e.g.
\cite{HMS07} for an example where this is not the case.
Furthermore, in the context of the $C^1$-stability theorem (see \cite{Ma78},
\cite{Ma86}), Ma\~n\'e conjectured that the $C^1$ interior of the
subset of diffeomorphisms whose periodic points are
hyperbolic is the set of diffeomorphisms satisfying Axiom A and the
 no-cycles condition. Such conjecture was proved by N.
 Aoki \cite{NA}. By Kupka-Smale theorem there exists a
 residual subset of $C^1$-diffeomorphisms exhibiting only
hyperbolic  periodic points. This shows that there is a large
 difference in assuming that a property of hyperbolicity in the periodic
 set  holds residually instead of requiring it for an open set
 in the space of $C^1$-diffeomorphisms.

We will say that a system $f$ has the {\em ergodic closing property}
if any $f$-ergodic measure can be weak* approximated by ergodic
measures supported in $f$-periodic orbits. By the Ergodic
Decomposition Theorem, this is obviously equivalent to say that
$\SM_1(f)$ is the closure of the convex hull of ergodic measures
supported in $f$-periodic orbits.

As a by-product of our techniques we also obtain the following
consequence:

\begin{corollary}{\label{cor1}}
Suppose that $f \in \DI(M)$ exhibits a NUH periodic set with
dominated splitting $T_{Per(f)}M= E^{cs} \oplus E^{cu}$. If  $f$ has
the {\em ergodic closing property} then, $\ov{Per(f)}$ is
hyperbolic.
\end{corollary}

In the remaining of this section we deal with the non-invertible setting.
As the notion of Axiom A endomorphisms is slightly different
and quite elaborated notion we deal first with the case of expanding transformations.
For completeness reasons let us recall some necessary definitions.

\begin{definition}
A $C^1$-map $g:M \to M$ on a compact manifold $M$
is {\em expanding} if there are constants $C> 0$ and $\sigma> 1$
such that
$$
\|[Dg^n(x)]^{-1}\|< C \cdot \sigma^{-n}, \forall n \in \natural.
$$
\end{definition}

The following notion deals with expansion at an asymptotic level and
it implies that all Lyapunov exponents of points in $X$ are positive.

\begin{definition}
We say that a map $g:M \to M$
is non uniformly
expanding ({\em NUE}) on a set $X \subset M$
if there exists $\la< 0$ such that
$$
\liminf\limits_{n \rightarrow \infty} \frac 1 n \sum_{j=0}^{n-1}\log
\|[Dg(g^j(x))]^{-1}\| \le \la<0\mbox{ for all }x\in X.
\label{defnue}
$$
\end{definition}

We are now in a position to state our criteria to
expansion for local diffeomorphisms.

\cmt \label{theo4} Let $\SS \subset \SU$ a residual subset of
an open set $\SU$ contained in $\NE(M)$.
Suppose that each $g \in \SS$ is non uniformly
expanding on the set $Per(g) \subset M$ of periodic points. Then,
there exists an open and dense subset of $\SU$ whose elements are
expanding maps.
\fmt

We also obtain a result analogous to Corollary \ref{cor1} above:

\begin{corollary}{\label{cor2}}
Suppose that $g \in \NE(M)$ exhibits a NUE periodic set.  If  $g$
has the {\em ergodic closing property} then $f$ is an expanding map.
\end{corollary}

Now we deal with the more subtle case of Axiom A endomorphisms and
prove a  result analogous to the one of Theorem~\ref{theo3}. One of
the main difficulties when dealing with non-invertible maps is the
possible existence of positively invariant sets (that is, sets which
are equal to their images) for an endomorphism that are not
negatively invariant. In fact, the single existence of an unstable
space is not uniquely determined: if the dimension of unstable space
$E^u(x)$ of a point $x$  is not the same of the ambient manifold
(which is the expanding map case),  for each choice of negative
branch in the pre-orbit of $x$ we may obtain a different unstable
space for $x$. This motivates to the use of unstable cone fields.

Let $X \subset M$ and let $E$ be  a subbundle of the tangent bundle
$TM$ restricted to $X$. Given a point $p \in X$, the cone $C_a(p)$
of width $0< a= a(p)< 1$ around $E(p)$ by
$$
C_a(p):= \{ v \in T_pM, \min_{w \in E(p)}\{\angle (v, w) \leq a\}\}
$$
We define the \emph{cone field $C_a$} of width $a:X \to (0, 1)$ around $E$
as the map $X \ni p \mapsto C_a(p)$. The cone field is continuous if
both $a$ and $p \mapsto E(p)$ are continuous.

We say that a cone field $C_a$ around $E$ and a subbundle $\hat E$
are complementary, if $E(p) \oplus \hat E(p)= T_pM$, and $C_a(p)
\cap \hat E(p)$ is trivial, for all $p$ in the intersection of the
domains of $E$ and $\hat E$.


We now present the definition of hyperbolic set for non-singular endomorphisms.

\begin{definition}
Let $\La$ be a positively invariant compact set for a $C^1$ nonsingular
endomorphism $g$ of a manifold $M$ and set $X:= \cup_{n=
0}^{+\infty} g^{-n}(\La)$. We say that $\La$ is a \emph{hyperbolic set} if
there are complementary invariant subbundle  $E^s$ of $T_\La M$ and
a positively invariant cone field $C^u_a$ defined on $T_XM$, such
that:
\begin{itemize}
\item $E^s$ is $Dg$-invariant,
that is,  $Dg(E^s(y))\subset E^s(g(y))$, $\forall y \in \La$;

\item $C^u_a$ is a $Dg$-invariant, ie, $Dg(x) \cdot C^u_a(x)
\subset C^u_a(g(x))$, $\forall x \in X$;

\item The angle between $E^s(y)$ and $C^u_a(y)$ is greater than a
positive constant, $\forall y \in \La$.

\item There are constants $c> 0$, $0< \vsi< 1$, such that $\forall n \in
\natural$
$$
\|Dg^n(y)|_{E^s}\| < c \cdot \vsi^n, \forall y \in \La \qquad
\text{and} \qquad \qquad \qquad \qquad \qquad \quad
$$
$$
 \qquad \|[Dg^{n}(x)|_{C^u_a(x)}
]^{-1}\|:= \inf_{ v \in C^u_a(x), \|v\|= 1} \{\|Dg^{n}(x) \cdot
v\|^{-1}\} < c \cdot \vsi^n, \forall x \in X.
$$

\end{itemize}
\end{definition}

We introduce the following definition of non uniformly hyperbolic set
for endomorphisms.

\begin{definition}
\label{defNUHend}  Let $g: M \to M$ be a nonsingular
endomorphism on a compact manifold $M$. We say that a positively
invariant set $\La \subset M$ is a {\em non uniformly hyperbolic
set} (or simply NUH set) if
\begin{enumerate}
\item
There  are $Dg$-invariant, complementary subbundle $E^{cs}$ whose
domain  is $\La$  and a cone field $C^{cu}$ on 
$X:= \cup_{n= 0}^{+\infty} g^{-n}(\La)$;

\item There exists $\la< 0$ and an
adapted Riemannian metric for which any point $p \in \La$  satisfies
$$
\liminf_{n \to +\infty}\frac{1}{n}\sum_{j= 0}^{n-1} \log
\|Dg(g^j(p))|_{E^{cs}(g^j(p))}\| \leq \la
$$
and
$$
\liminf_{n \to +\infty}\frac{1}{n}\sum_{j= 0}^{n- 1} \log
\|[Dg(g^{j}(p))|_{C^{cu}(g^{j}(p))}]^{-1}\| \leq \la
$$
\end{enumerate}
\end{definition}

Roughly, comparing to the invertible context this condition replaces
the existence of a uniquely determined expanding direction by means
of the existence of an unstable complementary cone field exhibiting
non-uniform expansion.

\cmt\label{theo5} Let $\SU \subset \NE(M)$ be an open subset of
nonsingular endomorphisms on a compact manifold $M$. Assume that for
any $g$ in some residual subset $\SS$ of $\SU$, the set $Per(g)$ of
periodic points of $g$ is non uniformly hyperbolic, that the
subbundle $E^{cs}$ extends continuously to $\ov{Per(g)}$ and that
the cone field $C^{cu}$ extends continuously to $\ov{\cup_{n=
0}^\infty g^{-n}(Per(g))}$. Then, there exists a residual subset
$\SA$ of $\SU$ whose elements are Axiom A endomorphisms. In
particular, $\SU$ is contained in the closure of the set of Axiom A
endomorphisms. \fmt

For the proof of the theorems in the endomorphisms case, we use

\cmt \label{theo6} {(Ergodic Closing Lemma for nonsingular
endomorphisms.)} There exists a
residual subset $\SR \subset \NE(M)$ such that for any $f \in \SR$,
the set of $f-$invariant probabilities $\SM_1(f)$ is the closed
convex hull of ergodic measures supported on periodic orbits of $f$.
\fmt

\section{Proof of the Criteria for generic Hyperbolicity}\label{s.p.criteria}

Throughout this section we assume that $f$ is a diffeomorphism in a $C^1$-open set $\SU$
and that the periodic set $Per(f)$ is non uniformly hyperbolic.

\begin{remark}\label{remper}
It is not hard to check 
that the set $Per(f)$ is non uniformly hyperbolic if and only if there exists $\la< 0$ such that
every periodic point $p$ with period $t(p)$ satisfies
$$
\sum_{j= 0}^{t(p)-1} \log(\|[Df|_{E^{cu}}(f^j(p))]^{-1}\|)\leq \la
\cdot {t(p)}
$$
and
$$
\sum_{j= 0}^{t(p)-1} \log(\|Df|_{E^{cs}}(f^j(p))\|) \leq \la \cdot
{t(p)}.
$$
\end{remark}

Before the  proof of Theorem~\ref{theo3} we recall some notations and
preliminary results. The first of these results is the $C^1$-Closing Lemma proved
by Pugh~\cite{P1}.

\begin{theorem}
There is a residual subset $\tilde\SR$ of $\DI(M)$ such that $\Omega(f) =
\ov{Per(f)}$, $\forall f \in \tilde\SR$.
\end{theorem}

Using this version of Pugh's Closing Lemma, Ma\~n\'e \cite{Ma82} stated the following
result, often known as Ergodic Closing Lemma,
which plays an important role in our proof.

\begin{theorem}
\label{teoma}
There is a residual $\SR \subset \DI(M)$
such that for each $f \in \SR$, the set $\SM_1(f)$ of $f$-invariant
probability measures is the closed convex hull of the ergodic measures supported
at hyperbolic periodic orbits of $f$.
\end{theorem}


Given an $f$-invariant probability measure $\mu$ it is well known
from the Oseledets Theorem \cite{Os68}  that the \emph{Lyapunov
exponent}  at $x$ in the direction $v\in T_x M$
$$
\lambda(x,v)= \lim\limits_{n\rightarrow \infty} \frac{1}{n} \log
\|Df^n(x)v\|,
$$
is well defined in a set of total probability that is, a set that of
full measure for every invariant probability measure. In \cite{YC},
Cao proved that every nonsingular endomorphism $g$ for which all
invariant measures have only positive Lyapunov exponents is an
uniformly expanding map. An analogous result for diffeomorphisms
admitting continuous splitting also follows.
We will use the following technical lemma:

\begin{lemma}\label{lema2} \cite[Lemma 2]{YC} Let $g:K \to K$ be a continuous map
defined in a compact metric space $K$.
Let $\:\mathcal{M}_1(g)$ be the space of $g$-invariant probabilities and
let $\phi$ be
a continuous function on $K$. If $\int{\phi d\mu}<{\lambda},\forall
\ \mu \in \:\mathcal{M}_1(g)$, then  there exists $N \in \natural$
such that $\forall n\geq N$ we have
$$
\frac 1 n \sum _{i=0}^{n-1} \phi(g^i(x)) < \lambda, \forall x \in K.
$$
\end{lemma}

Notice that if $g \in \NE(M)$ and the condition in hypothesis
of the lemma holds for $\phi=\log\|Dg^{-1}\|$ then $g$ is expanding.
Analogously, if $f \in \DI(M)$, $\Lambda$ is a compact invariant set
with the continuous splitting $T_\La M=E^{cs}\oplus E^{cu}$,
and the functions $\phi_1= \log\|Df|_{E^{cs}}\|$, $\phi_2=\log\|Df^{-1}|_{E^{cu}}\|$
satisfy the condition of the lemma then $\La$ is a hyperbolic set for $f$.  Indeed, this
is an immediate  consequence of the following simple proposition:

\begin{proposition} \label{propr}
Let $g: K \to K$ be a continuous map defined in a compact metric
space $K$. Let $\phi_n: K \to \re$ be a $g$-subadditive sequence of
continuous functions (that is, $\phi_{n_1+ n_2}(x) \leq
\phi_{n_1}(x)+ \phi_{n_2}\circ g^{n_1}(x)$, $\forall n_1, n_2 \in
\natural$). Suppose that $\int_M \phi_1 d\mu < \la$, $\forall \mu
\in \SM_1(g)$. Then there exists $N$ such that for all $n \geq N$,
$$
\phi_n(x) \leq n \cdot \la, \forall x \in K.
$$
\end{proposition}

\begin{proof}
By lemma \ref{lema2} applied to $\phi= \phi_1$, there is $N \in
\natural$ such that
$$
\frac 1 n \sum_{i=0}^{n-1} \phi_1(g^i(x)) < \lambda, \forall x \in
K.
$$
The subadditivity of $\phi_n$ then implies that
$$
\phi_n \leq \sum_{i=0}^{n-1} \phi_1(g^i(x)) < n \cdot \la, \forall n
\geq N, \forall x \in K.
$$
\end{proof}

Taking into account these preliminar results we proceed to prove
the main results stated in the previous section.

\begin{proof}[Proof of Theorem~\ref{theo3}]
Let $\tilde\SR$ and $\SR$ be the residual sets given by
Pugh's and Ma\~n\'e's Closing Lemmas, respectively, and let
$f:M \to M$ be a diffeomorphism in the residual set
$\SA$ given by  the intersection of $\SS$, $\tilde\SR$ and $\SR$.
Therefore, not only $\Omega(f)= \ov{Per(f)}$ as
 the space of invariant probabilities for $f$ is the closed
convex hull of ergodic probabilities supported in periodic orbits.
Take $\phi_s:= \log(\|Df|_{E^{cs}}\|)$, $\phi_n= \log(\|Df^n|_{E^{cs}}\|)$,
$\phi_u:=\log(\|[Df|_{E^{cu}}]^{-1}\|)$, $\psi_n= \log(\|[Df^n|_{E^{cu}}]^{-1}\|)$
for $n\in\mathbb N$ and $K= \Omega(f)$ in Proposition~\ref{propr}. Observe that such functions
are continuous since the subbundles $E^{cs}$ and $E^{cu}$ vary continuously.

Since $f \in \SS$, we can use Remark~\ref{remper} to deduce that there exists $\la< 0$
such that for every periodic point $p$ with period $t(p)$ we have
$$
  \int_{\Omega(f)}
\phi_s d \big(_{\frac{1}{t(p)} \sum_{j= 0}^{t(p)-
1}\delta_{f^{j}(p)}}\big)= \frac{1}{t(p)} \sum_{j= 0}^{t(p)- 1}
\phi_s(f^j(p)) \leq \la
$$
and
$$
\int_{\Omega(f)} \phi_u d \big(_{\frac{1}{t(p)} \sum_{j= 0}^{t(p)-
1}\delta_{f^{j}(p)}}\big)= \frac{1}{t(p)} \sum_{j= 0}^{t(p)- 1}
\phi_u(f^j(p)) \leq \la,
$$
where $\frac{1}{t(p)}\sum_{j= 0}^{t(p)- 1} \delta_{f^j(p)}$ is the
ergodic measure supported in the periodic orbit of $p$. Using the Ergodic
Closing Lemma for $f$, all $f$-invariant
probability $\mu$ is the limit of a convex combination of such
measures supported in periodic orbits and so we conclude that
$$
\int_{\Omega(f)} \phi_s d\mu \leq \la
\quad\text{ and }\quad
\int_{\Omega(f)} \phi_u d\mu \leq \la, \forall \mu \in \SM_1(f).
$$
Note that $\SM_1(f|_{\Omega(f)})\simeq \SM_1(f)$.
We are under the hypothesis of Proposition   \ref{propr} (just
exchange $\la$ by $0> \la'> \la$ to guarantee the strict inequality in the
statement of lemma \ref{lema2}). Hence, we conclude that
$\ov{Per(f)}= \Omega(f)$ is a hyperbolic set for $f$.

Now, let us see that (residually) $\Omega(f)$ has no cycles, and therefore, by \cite{Pa87},
$f$ is $\Omega-$stable. In particular, there is a neighborhood of $f$ whose elements are
Axiom A.

Given a set $S \subset M$ let $S(\epsilon)$ be the $\epsilon-$neighborhood of $S$ in $M$.

Consider $\SC$ the collection of compact subsets of $M$ endowed with Hausdorff distance given by
$$
d_H(K_1, K_2):= \inf\{\epsilon> 0, K_1 \subset K_2(\epsilon) \text{ and } K_2 \subset K_1(\epsilon)\},
\forall K_1, K_2 \in \SC,
$$

Let $\SK$ be the set of Kupka-Smale diffeomorphisms, which is residual.
Since hyperbolic periodic points are robust, the map $\Psi: \SK \to \SC$ that assigns each $f$ to $\ov{Per(f)}$ is lower semicontinuous.
Hence, there exists a residual subset $\SR' \subset \SU$ of continuity
points for $\Psi$.

Note that if $g \DI(M)$ is an Axiom A with cycles, then $g$ is a discontinuity point for
$\Psi$.  In fact, by \cite{Pa70}, there exists an wandering point $y$ out from a neighborhood of $\Omega(g)$
which by an  arbitrarily small perturbation can be changed into
a point of transversal intersection between Stable and Unstable Manifolds of some basic set
of $g$. Such point $y$  then becomes a nonwandering point for $g_n$ in a sequence
tending to $g$ as $n \to +\infty$, and is contained in $\ov{Per(g_n)}$.
Therefore, $\ov{Per(g_n)} \not\to \ov{Per(g)}$ as $n \to +\infty$.
Since $\SK$ is a residual subset of $\DI(M)$, and transversal intersection between submanifolds
is an open property, we can take $g_n$ to be Kupka-Smale diffeomorphisms.
This implies that $g$ is not a point of continuity for $\Psi$.

So we conclude that for $f$ in a residual subset of $\SU$,
$f$ is Axiom A and has no cycles. Therefore, we obtained an open and dense
subset  $\SA \subset \SU$,  whose elements are Axiom A (with no cycles).

\end{proof}

\begin{remark}
Under the conditions of either Theorem~\ref{theo2} or
Theorem~\ref{theo3}, if it occurs that for $f \in \SA$ is such that
$\Omega(f)$ has the no cycle condition (see \cite{Pa70}) then we
conclude that $\Omega(\hat f)$ is hyperbolic for all $\hat f$
belonging in an open and dense subset of $\SU$.
\end{remark}

As we discussed in the previous section the results in Theorem~\ref{theo2}
follow from  Theorem~\ref{theo3} by proving that the existence of a dominated splitting
over the periodic set is continuous.  More precisely, it is enough to recall the following lemma.

\begin{lemma}\label{ledom} \cite[Lemma 14]{COP}
Let $f:M \to M$ be a diffeomorphism on a compact manifold $M$. Let
$X\subset M$ be some $f-$invariant set. Suppose there exists some
invariant dominated splitting $T_X M= E \oplus \hat E$. Then, such
splitting is continuous in $T_X M$, and unique once we fix the
dimensions of $E, \hat E$. Moreover, it extends uniquely and
continuously to a splitting of $T_{\ov X}M$.
\end{lemma}

Note that Corollary \ref{cor1} is an immediate consequence of the
proofs of Theorems \ref{theo3} and \ref{theo2} above, as in such
corollary we are assuming that the thesis in Ergodic Closing Lemma
holds for $f \in \DI(M)$.

\begin{remark}\label{rempernue}
Just as in the case of Remark~\ref{remper},  $Per(g)$ is NUE if and only if there
exists $\la< 0$ so that  $\sum_{j= 0}^{t(p)-1} \log(\|[Dg(g^j(p))]^{-1}\|)\leq
\la \cdot {t(p)}$ for each periodic point $p$ with period
$t(p)$.
\end{remark}

The remainder of this section is now devoted to the proof of Theorem~\ref{theo4} and \ref{theo5} in the non-invertible
context. As the reader may guess the proof of these theorems goes
along the same lines and ideas used in the proofs of Theorems
\ref{theo2} and \ref{theo3},
 using the Ergodic Closing Lemma version for Nonsingular Endomorphisms (Th. \ref{theo6}).
 For completeness we write down the proofs anyway.

\begin{proof}[Proof of Theorem~\ref{theo4}]
Let $g:M \to M$, $g \in \SS \subset \NE(M)$ belonging in the
residual set $\SR$ given by Th. \ref{theo6}.
Therefore, in particular, we have that the space of
invariant probabilities for $g$ is the closed convex hull of ergodic
probabilities supported in periodic orbits.  Put $\phi:=
\log(\|[Dg]^{-1}\|)$. Since $g \in \SS$, there exists $\la< 0$ such
that, for any periodic point $p$ with period $t(p)$ we have (see
Remark \ref{rempernue}):
$$
\int_M \phi d \big(\frac{1}{t(p)} \sum_{j= 0}^{t(p)-
1}\delta_{g^{j}(p)}\big)= \sum_{j=0}^{t(p)- 1} \phi(g^j(p))\leq \la,
$$
where $\frac{1}{t(p)} \sum_{j= 0}^{t(p)- 1} \delta_{g^j(p)}$ is the
ergodic measure supported in the periodic orbit of $p$. As the
thesis of Ergodic Closing Lemma holds for $g$, all $g$-invariant
probability $\mu$ is the limit of a convex combination of such
measures supported in periodic orbits and so we conclude that
$$
\int_M \phi d\mu \leq  \la, \forall \mu \in \SM_1(g).
$$
We are again under the hypothesis of the fundamental proposition
\ref{propr}, which implies that $g$ is an expanding map.

\end{proof}

Corollary \ref{cor2} is a straightforward consequence of the proof
above.

\begin{proof}[Proof of Theorem~\ref{theo5}]

Let $g:M \to M$ be a nonsingular endomorphism belonging to the
residual set $\SA$ given by  the intersection of $\SS$, $\tilde\SR$
and $\SR$, the last ones given by the  Closing Lemmas for
Endomorphisms. Therefore, in particular, we have $\Omega(g)=
\ov{Per(g)}$ and the space of invariant measures for $g$ is the
closed convex hull of ergodic measures supported in periodic orbits.
Put $\phi_s:= \log(\|Dg|_{E^{cs}}\|)$, $\phi_u:=
\log(\|[Dg|_{C^{cu}}]^{-1}\|)$. Such functions are continuous, since
the subbundle $E^{cs}$ and the cone field $C^{cu}$ are assumed to be
continuous. $\phi_s$ and $\phi_u$ are defined in $\Omega(g)$ and $X=
\ov{\cup_{n \geq 0} g^{-n}(Per(g))}$, respectively. Since $g \in
\SS$, there exists $\la< 0$ such that, for any periodic point $p$
with period $t(p)$ we have (just as in Remark \ref{remper}):
$$
  \int_{\Omega(g)}
\phi_s d \big(_{\frac{1}{t(p)} \sum_{j= 0}^{t(p)-
1}\delta_{g^{j}(p)}}\big)= \frac{1}{t(p)} \sum_{j= 0}^{t(p)- 1}
\phi_s(g^j(p)) \leq \la
$$
and
$$
\int_X \phi_u d \big(_{\frac{1}{t(p)} \sum_{j= 0}^{t(p)-
1}\delta_{g^{j}(p)}}\big)= \frac{1}{t(p)} \sum_{j= 0}^{t(p)- 1}
\phi_u(g^j(p)) \leq \la,
$$
where $\frac{1}{t(p)}\sum_{j= 0}^{t(p)- 1} \delta_{g^j(p)}$ is the
ergodic measure supported in the periodic orbit of $p$. The Ergodic
Closing Lemma  for Nonsingular Endomorphisms implies that
$$
\int_{\Omega(g)} \phi_s d\mu \leq \la \quad  \text{ and } \quad
\int_X \phi_u d\mu \leq \la, \forall \mu \in \SM_1(g).
$$
Just as in the previous theorems, we apply Proposition   \ref{propr}
 to conclude that $\ov{Per(g)}= \Omega(g)$ is a hyperbolic set for
$g$.

\end{proof}

\section{Proof of the Ergodic Closing Lemma for nonsingular endomorphisms}\label{s.ecl}

\label{sec3} For the proof of Theorem \ref{theo6},  let us start
by fixing some notation. Given $x \in M$,
we define $B_\epsilon(f, x)$ as an $\epsilon-$neighborhood of the
orbit of $x$. Define  $\Sigma(f)$ as the set of points $x \in M$
such that for every neighborhood $\SU$ of $f$ and every $\epsilon>
0$, there exist $g \in \SU$ and $y \in M$ such that $y \in Per(g)$,
$g= f$ on $M \setminus B_\epsilon(f, x)$ and $d(f^j(x), g^j(y)) \leq
\epsilon$, $\forall 0 \leq j \leq m$, where $m$ is the $g-$period of
$y$.

At the end of this section, we shall obtain the {\bf residual
} version of the Ergodic Closing Lemma (our Theorem \ref{theo6}) as consequence of
the following result:
\begin{theorem}\label{theop}
For any nonsingular endomorphism $f$, $\Sigma(f)$ is a total probability set,
that is, $\Sigma(f)$ is a full probability set for any $f-$invariant
probability.
\end{theorem}

\begin{remark}
Theorem \ref{theop} above, with $f \in \DI(M)$ instead of $f \in
\NE(M)$ in its statement, was the former Ergodic Closing Lemma
proved by Ma\~n\'e. In fact, Ma\~n\'e, in \cite{Ma82}, did not
explicitly give the proof of his corresponding classical residual
version, although he stated such version in \cite{Ma87}. This gap was
filled out very recently by  Abdenur et al in \cite{ABC}.
\end{remark}

\begin{definition}{($\epsilon-$shadowing by a periodic point.)} Let $f$ and $g$ maps on
a compact metric space $\La$. Given $\epsilon> 0$ and $x \in \La$, we say that
a $g-$periodic point $p$ with period $n$  {\em $\epsilon-$shadows $x$} iff
$d(g^j(p), f^j(x))< \epsilon, \forall 0 \leq j \leq n$.
\end{definition}

Let $\epsilon> 0$ and a neighborhood $\SU \ni f$, $\SU \subset \NE(M)$ be given.
 We define {\em $\Sigma(f, \SU, \epsilon)$} as
the set of points $x \in M$ such that there exist $g \in \SU$ and $y \in
M$ such that $y \in Per(g)$, $g= f$ on $M \setminus B_\epsilon(f,
x)$ and $d(f^j(x), g^j(y)) \leq \epsilon$, $\forall 0 \leq j \leq
m$, where $m$ is the $g-$period of $y$. That is, $\Sigma(f, \SU, \epsilon)$ is
the set of points $x \in M$ which are $\epsilon-$shadowed by a
periodic point $y \in Per(g)$, for some $g \in \SU$.
Everytime there is no chance of misunderstanding, we will just
write $\Sigma(\SU, \epsilon)$ instead of $\Sigma(f, \SU, \epsilon)$.
If we take a nested neighborhood basis $\SU_n$ of $f$ in $\NE(M)$ then
$$
\Sigma(f)= \cap_{n \in \natural} \Sigma(f, \SU_n, 1/n).
$$
Therefore, Theorem \ref{theop} is an immediate consequence of

\begin{proposition}\label{prop}
For any nonsingular endomorphism $f$, any neighborhood $\SU$ of $f$
and $\epsilon> 0$, $\Sigma(\SU, \epsilon)= \Sigma(f, \SU,\epsilon)$
is a total probability set for $f$.
\end{proposition}

Since the proof of Proposition \ref{prop} is long and involving
we divide it in two main parts that we shall describe next
for the  reader's convenience.
Part 1 consists in an improvement of Closing Lemma (see \cite{P1},
\cite{P2}, \cite{PR}, \cite{LW1}), roughly stating in particular that given a nonsingular
endomorphism
 $f$, for any point $x \in M$ there are $\rho> 1$ and $r> 0$ such that if $x$  returns to a  ball $B(x, r')$, $r'< r$, then
there is an intermediate iterate $y= f^{m(x)}(x)$ which returns to $B(x, \rho r')$ that is shadowed by a $g$-periodic point, where $g$ is some endomorphism close to $f$.
The precise statement corresponds to Lemma
\ref{lesomb}, whose proof we write down further in this paper. Let us point
out that this first part of the argument is entirely topological.

Part 2 uses ergodicity and Birkhoff's Theorem to show that the set of points $x \in M$ that
are shadowed by periodic points of some endomorphism close to $f$ has total probability with respect to $f-$invariant measures. In fact, recall first that typical points for an ergodic invariant measure are recurrent.
The arguments in Part 1 yield that at every time a recurrent point $x$ returns to a ball there is an intermediate iterate $y$ with the shadowing property we look for, that also returns to a ball with the same center and whose radius is
$\rho$ times bigger. Notice that if we could take $\rho$ equal to one above there would be nothing to do: the set of
recurrent points and of the $\epsilon-$shadowable points $\Sigma(\SU, \epsilon)$ would visit
equally the same balls which imply that they have the same measure.
More generally, using standard Radon-Nikodym derivative calculations, the same conclusion would be obtained if the $\mu-$measure of a ball and
its $\rho-$homothetic image were proportional, at least for balls with sufficiently small radius.
In fact, given a set $S$ it is well known that $\mu-$a.e. point $x \in M$ we have
$$
\chi_S(x)= \lim_{r \to 0} \frac{\mu(B(x, r) \cap S)}{\mu(B(x, r))}.
$$
Applying Part 1 when $S= \Sigma(\SU, \epsilon)$ one has that for some $\rho= \rho(x)$
$$
\liminf_{r \to 0} \frac{\mu(B(x, \rho r) \cap S)}{\mu(B(x, r))} \geq 1
$$
holds $\mu-$a.e. $x \in M$,  due to ergodicity of $\mu$ and the fact that asymptotic frequency of visits of a typical point
to $B(x, \rho r) \cap S$ is greater or equal than to $B(x, r)$.
If the measures of $B(x, r)$ and $B(x, \rho r)$ were proportional, say, by a factor of $\zeta> 0$,
we could obtain
$$
\chi_S(x) \geq \liminf_{r \to 0} \frac{\mu(B(x, \rho r) \cap S)}{\mu(B(x, \rho r))} \cdot \frac{\mu(B(x,  r))}{\mu(B(x,  r))} \geq \zeta> 0,
$$
which implies that the indicator function $\chi_S$ of $S$ is equal to $1$ a.e., and so
$S$ has total probability.
Even though such proportionality on balls with same center and different radius
may not hold for general measures,  rough estimates in such direction do hold, which are sufficient
for our proof. Following the ideas above, we use a  Vitali's covering
like argument to prove that the set of points that are shadowed by a
periodic point of some nearby endomorphism $g$ has total probability
for $f$. The core of these arguments is contained in Lemma
\ref{lerg} and Theorem \ref{teou} below.

We proceed to prove Part 1 and its main technical Lemma~\ref{lesomb}.
For that purpose we shall first introduce some notations and proceed to some
perturbation lemmas. As the manifold $M$ is compact, there is $\delta$ such
that $\{\exp_p, p \in M\}$ is an equilipschitz family of
diffeomorphisms, such that each exponential map $\exp_p$ embeds
$B(0, \delta)$ in a neighborhood $B_p$ of $p$. Given $p \in M$, we
define a metric $d'= d'_p:B_p \times B_p \to [0, +\infty)$ given by
$$
d'(x, y):= |\exp^{-1}_p(x)- \exp^{-1}_p(y)|.
$$

 Obviously, $d'$ is Lipschitz-equivalent to the manifold usual
metric restricted to $B_p$. Setting $d'$ as the metric in $B_p$,
then $\exp_p$ isometrically maps $B(0, \delta)$ on $B_p= B'(p,
\delta)$, where the quote ' signs the ball in the metric $d'$.

\begin{lemma}{\cite{PR}}\label{lelift}
For any $\eta > 0$, for any $f \in \NE(M)$, there is an $\alpha > 0$ such that for any $f
\in \NE(M)$, any $q \in M$, any two points $v_1, v_2 \in T_qM$ with
$B(v_2, |v_1- v_2|/\alpha) \subset B(0, \delta) \subset T_qM$, there
is a diffeomorphism $h= h_{q, \alpha, v_1, v_2}: M \to M$, called an
$\alpha-$lift, such that:
\begin{enumerate}
\item $h(\exp_q(v_2))= \exp_q(v_1)$;
\item The closure of set of points where $h$ differs to the identity
is contained in $\exp_q(B(v_2, |v_1- v_2|/\alpha)$;
\item $d_1(h f, f)< \eta$.
\end{enumerate}
\end{lemma}

\begin{definition}{(Dynamical neighborhood.)}
We say that a neighborhood $V$ of a point $p \in M$ is  $N-${
dynamical} for $f$ if each connected component $\cup_{j= 0}^N
f^{-j}(V)$ contains exactly one point of $\cup_{j= 0}^N
f^{-j}(\{p\})$.
\end{definition}

\begin{lemma}{\cite{F}, \cite{LW1}}\label{lelinea}
Let $f \in \NE(M)$, $p \in M$, $N \geq 1$ given such that all terms in
$\cup_{j= 0}^{N+ 1} f^{-j}(p)$ are distinct. Then, for any $\eta>
0$, there is a $\beta> 0$, and a map $f_1 \in \NE(M)$, called a local
linearization of $f$ with the following properties (1)-(5).
\begin{enumerate}
\item
$B'(p, \beta)$ is  $(N+ 1)$-{\em dynamical } for both $f$ and $f_1$,
and $f^{-j}(B'(p, \beta))= f_1^{-j}(B'(p, \beta))$ for $j= 1, \dots,
N+1$.

\item For $q \in \cup_{j= 1}^{N+ 1} f^{-j}(p)$, let $V(q)$
be the open connected component of $\cup_{j= 1}^{N+ 1} f^{-j}(B'(p,
\beta/4))$ containing $q$. Then, $f_1|_{V(q)}= \exp_{f(q)} \circ
(T_q f) \circ \exp_q^{-1}$.

\item $f_1^{N+1}(x)= f^{N+1}(x)$, $\forall x \in f^{-N-1}(B'(p,
\beta))$.

\item $f_1= f$ on $M \setminus \cup_{j= 1}^{N+ 1} f^{-j}(B'(p,
\beta))$.

\item $d_1(f_1, f)< \eta$.

\end{enumerate}

\end{lemma}

\begin{remark} \label{rent}
For the sequel, we need to emphasize two aspects from the last lemma.
On the one hand, we obtain as a direct consequence that if for some $k
\in \natural$ we have that $x, f^k(x)$ are both out of $\cup_{j=
0}^{N} f^{-j}(B'(p, \beta))$, then $f^k(x)= f_1^k(x)$. On the other hand,
notice that as $\eta$ goes to $0$, one can take $\beta$ arbitrarily small in
the last lemma.
\end{remark}

\begin{theorem}{(Theorem A in  \cite{LW1}.)} \label{teoalw}
Let $(\ST, T_q)$ a complete tree of isomorphisms associated to the
pre-orbit of a point $q_0 \in M$, that is, a collection of
$n-$dimensional inner product spaces $E_{q}$ and isomorphisms $T_q:
E_q \to E_{q_0}$ associated to each $q$ in the pre-orbit of $q_0$,
with $T_{q_0}$ equal to identity. Given $\alpha> 0$, there are
$\rho> 0$ and $N \geq 1$ such that: for any ordered set $X= \{x_0
\prec \dots \prec x_t \} \subset E_{q_0}$, there is a point $y \in X
\cap B(x_t, \rho |x_0- x_t|)$ such that for any branch $\Ga= \{q_0,
q_1, \dots, \}$ of $\ST$, there is a point $w= w(\Ga) \in X \cap
B(x_t, \rho |x_0- x_t|)$ which is before $y$ in the order of $X$,
together with $N+ 1$ points $c_0(\Ga)= c_0, \dots, c_N(\Ga)= c_N \in
B(x_t, \rho |x_0- x_t|)$ satisfying the following two conditions:
\begin{itemize}
\item $c_0= w, c_N= y$; and

\item $|T_{q_n}^{-1}(c_j)- T_{q_n}^{-1}(c_{j+1})| \leq \alpha
d(T_{q_n}^{-1}(c_{j+1})- T_{q_n}^{-1}(A))$, where $A:= \{x \in X,
w\prec x \prec y\} \cup \partial B(x_t, \rho |x_0- x_t|)$.

\end{itemize}
\end{theorem}

The next lemma is the main target in the first part of our
Ergodic Closing Lemma, whose arguments are just topological.
It implies in particular that, given $\epsilon> 0$ and
  {\em any} $f-$recurrent
point $x$, then $x$ has an iterate which is $\epsilon-$shadowed by a periodic point of some
$g$ close to $f$.

\begin{lemma}
\label{lesomb} Given $f \in \NE(M)$, $p \in M$, $\epsilon> 0$ and a
neighborhood $\SU$ of $f$, there exist $r> 0$, $\rho'> 1$ such that
if for some natural $t> 0$, we have $x, f^t(x) \in B'_{\ov r}(p)$,
with $0< \ov r \leq r$, then there exist $0 \leq t_1 < t_2 \leq t$
and $g \in \SU$ such that:
\begin{itemize}
 \item $w= f^{t_1}(x), y= f^{t_2}(x) \in \ov{B'_{\rho' \ov r}(p)}$;
 \item $g^{t_2- t_1}(w)= w$;
 \item $g(z)= f(z)$ for $z \notin B_\epsilon(f, x)$ and $d(g^j(w), f^j(w)) \leq
\epsilon$, $\forall 0 \leq j \leq t_2- t_1$.
\end{itemize}
\end{lemma}

\begin{proof} Take an $\eta> 0$ such that the $\eta-$ball with center $f$ is
contained in $\SU$.
 Take $1> \alpha> 0$ such that $d_1(h \circ f, f) < \eta/2$, for any
$\alpha-$lift $h$.
Without loss of generality, we assume that $\epsilon< \alpha^2$. We
also assume that $\epsilon< \delta$.

We assume that $p$ is not periodic for $f$, otherwise, there is
nothing to prove.
This implies that all points in the pre-orbit of $p$ are distinct.
Let $\rho> 2$ and $N \geq 1$ be the numbers provided  by Theorem
\ref{teoalw}, for $\alpha> 0$ taken as above, and for $q_0= p$, each
$q_j$ to be some $j-$pre-image of $p$,  $E_{q_j}= T_{q_j} M$ and
$T_{q_j}= Df^j(q_j)$.

So, take $r> 0$ such that $r< \epsilon/(6\rho)$ and
$\diam(f^{-j}(B'(p, 3\rho r))< \epsilon$, $\forall j= 0, \dots,
N+1$. We assume that each connected component of $\cup_{j= 0}^{N+ 1}
f^{-j}(B'(p, 3\rho r))$ contains exactly one point $q_j \in
f^{-j}(p), j= 1, \dots, N+1$. In particular, if $z, f^{\hat t}(z) \in
B'(p, 3\rho r)$, then $\hat t> N+ 1$.

Now, assuming that $x, f^t(x) \in B'(p, \ov r)$, for some $0< \ov r
< r$ we can apply Theorem \ref{teoalw} to the set $X= \{x, f(x),\dots,
f^t(x)\} \cap B'(p, 3\rho \ov r)$ endowed with the order given by
the iterate number: if $f^k(x), f^{\hat k}(x) \in X$, then $f^k(x)
\prec f^{\hat k}(x) \Leftrightarrow k< \hat k$. Therefore, set
$\rho'= 3 \rho$. We then obtain $f^{t_2}(x)= y \in \{x, \dots,
f^t(x)\} \cap B'(f^t(x), \rho \cdot d'(f^t(x), x)) \subset B'(p,
\rho' \ov r)$ such that for any branch $\Ga= \{p= p_0, p_1, \dots,
p_n, \dots\}$ of the pre-orbit of $p$, there is $w= w(\Ga)=
f^{t_1}(x) \in \{x, \dots, f^t(x)\} \cap B'(f^t(x), \rho \cdot
d'(f^t(x), x))$, with $t_1= t_1(\Ga)< t_2$ together with points
$c_0= c_0(\Ga), \dots, c_N= c_N(\Ga) \in B'(f^t(x), \rho \cdot
d'(f^t(x), x))$ such that:

\vspace{0.1cm}

{(a)} $c_0= w, c_N= y$; and

{(b)} $|T_{p_j}^{-1}(c_j)- T_{p_j}^{-1}(c_{j+1})| \leq \alpha
d(T_{p_j}^{-1}(c_{j+1}), T_{p_j}^{-1}(A))$, where $A:= \{f^j(x) \in
X; t_1< j< t_2 \} \cup \partial B'(f^t(x), \rho \cdot d'(f^t(x),
x))$.

\vspace{0.2cm}

As $w= w(\Ga)$ and $y$ are both in $X$, there is a natural number
$k(\Ga) \geq 1$ such that $f^{k(\Ga)}(w(\Ga))= y$. Note that
$k(\Ga)> N+ 1$, as $\cup_{j= 0}^{N+1} f^{-j}(\{y\}) \cap B'(p, 3\rho
r)= y$, from our choice of $r$. Setting $z:= f^{k(\Ga)- N-
1}(w(\Ga))$, we see that $z$ does not depend on the branch $\Ga$ of
$p$, since $w(\Ga)$ and $y$ are in $X$, $f^{k(\Ga)}(w(\Ga))= y$ and
$y, N$ do not depend on $\Ga$. By our choice of $r$, since $y \in
B'(p, 3\rho r)$ there is a unique connected component $V_{N+1}
\subset f^{-(N+1)}(B'(p, 3\rho r))$ such that $z \in V_{N+1}$. Also,
there is a unique $p_{N+1} \in f^{-(N+1)}(p) \cap V_{N+1}$. From now
on, we fix $\Ga$ as some branch of $p$ containing $p_{N+1}$ (That
is, $\Ga= (p= p_0, \dots, p_{N+1}, \dots)$), and we consider all
constants $w, c_0, \dots, c_N$, $k$ obtained by applying Theorem
\ref{teoalw} with respect to such branch. For each $p_j, j= 0,
\dots, N-1$, let $h_{p_j}$ be the $\alpha-$kernel lift obtained by
treating in lemma \ref{lelift} $q= p_j$, $v_1=
[Df^j(p_j)]^{-1}(c_j)$, $v_2= [Df^j(p_j)]^{-1}(c_{j+1})$. Defining a
map $g: M \to M$ by
$$
g:= \left\{
\begin{matrix}
h_{p_j} \circ f_1 \text{ on } V(p_{j+1}); \\
f_1    \text{ on the rest of } M,
\end{matrix}
\right.
$$
we have that $g \in \NE(M)$ and $d_1(g, f)< \eta$. Thus $g \in \SU$.

Due to condition (b) above, the $g-$orbit from $w$ to $z$ never
touches the region in which $g \neq f_1$. Therefore, $g^{k-
(N+1)}(w)= f_1^{k- (N+1)}(w)$. By remark \ref{rent}, we also have
that $f^{k- (N+1)}(w)= f_1^{k- (N+1)}(w)$, and thus
$$
g^{k- (N+1)}(w)= f^{k- (N+1)}(w)= z.
$$
Now, it is easy to see that $g^{N+1}(z)= w$ and then $g^k(w)= w$. In
fact, $f_1^{N+1}(z)= y$, and the lifts $h_{p_{N-1}}, \dots h_0$
gradually and slightly modifies $f_1-$orbit segment joining $z$ and
$y$, in such way that $g^{N+1}(z)= w$ and $d(g^j(z), f^j(z)) \leq
d(g^j(z), f_1^j(z))+ d(f_1^j(z), f^j(z))< \epsilon$, $\forall j= 1,
\dots, N+1$.

\end{proof}

Now we proceed to prove the second part of the proof of Proposition \ref{prop}.
Although   the main idea of this part is borrowed from \cite{Ma82},
the proofs we have written are presented in an abstract setting for
future use and bookkeeping purposes. This will also clarify the sort
of arguments which are used. 

We start by introducing some notation. We say that a subset $C$ of the torus $T^s$ is a cube
if it can be written as $A= I_1 \times \dots \times I_s$, where the sets $I_i$ are
intervals of same length in $S^1$ (containing both, none, or one of its boundary points).
If $p_i$ is the middle point of $I_i$, we say that the point $(p_1, \dots, p_s)$
is the center of $A$. The length of the intervals $I_i$ is called the side of the cube.
For each $k \in \natural^+$, let $(\SP^{(k)}_j)_{j \in \natural^+}$ be a  sequence of partitions
of $T^s$ by cubes whose side is $2\pi/k^j$. For every atom $P$ of a partition $\SP^{(k)}_j$, we can associate
cubes $\hat P$ and $\tilde P$ having the same center of $P$, but with sides $2\pi/k^{j-1}$
and $6\pi/k^{j-1}$, respectively. If $x \in T^s$, denote by $P^{(k)}_j(x)$
the atom of $\SP^{(k)}_j$ containing $x$. Suppose that $M$ is isometrically embedded in $T^s$.
We recall the following useful fact on such kind of partitions.

\begin{lemma}\label{leprp}\cite[Lemma~I.5]{Ma82}
For every probability measure $\mu$ on the Borel sets of $T^s$, every $\delta> 0$ and for all
odd natural $k$, the following inequalities holds for any $j \geq 1$:
$$
\mu(\{x; \mu(P^{(k)}_j(x)) \geq \delta \mu(\hat P^{(k)}_j(x))\}) \geq 1- \delta k^s
$$
and
$$
\mu(\{x; \mu(P^{(k)}_j(x)) \geq \delta \mu(\tilde P^{(k)}_j(x))\}) \geq 1- \delta 3^sk^s.
$$
\end{lemma}

Let $f \in \NE(M)$, $\epsilon> 0$, a neighborhood $\SU$ of $f$
and an ergodic $\mu \in \SM_1(f)$ be given. Extend $\mu$ to a measure on $T^s$ by $\mu(A):= \mu(A \cap M)$,
for all Borel set $A \subset T^s$.
Let $\mho \subset M$ be some Borelian set and suppose that  $\mho(
r, \rho)$, where $r> 0$, $\rho> 1$, is some Borelian set whose
elements are points $x \in M$ with the following property:  if $y
\in B'_{r'}(x)$ for some $0< r' \leq r$ and $f^t(y) \in B'_{r'}(x)$,
for some $t> 0$ then there exist $0 \leq t_1< t$, such that
$f^{t_1}(y) \in \ov{B'_{\rho r'}(x)} \cap \mho$. { Take $r_i> 0$,
$\rho_l > 1$ two monotone sequences} converging respectively to 0
and $+\infty$.

Our first target in this second part  is to obtain an abstract
result (Theorem \ref{teou})  which will be essential in both  proofs of
Proposition \ref{prop} and Theorem \ref{theo6}.
Such result says that, if $\cup_{i, l} \mho(r_i, \rho_l)= M$, then
$\mho$ has total probability for $f$.

\begin{remark}
All results from this point of the paper up to Theorem \ref{teou} do not use much regularity of $f$. In fact, specifically for the statements from Lemma \ref{lerg} up to Theorem \ref{teou}, we only request  $f:M \to M$ to be
a Borelian map such that $\SM_1(f) \neq \emptyset$. This occurs, for instance, if $f$
is a continuous map.
\end{remark}

For each pair $(i, l)$, we can find
and odd natural $k= k(i, l)$ and $j(i, l)$ such that $\forall j \geq j(i, l)$ and $x \in T^s$ there exists $0 \leq r \leq r_i$ satisfying
$$
P^{(k)}_j(x) \subset B_r(x)
$$
and
$$
\hat P^{(k)}_j(x) \supset \ov{B_{\rho_l r}(x)},
$$
the balls here are taken in the torus.
The next lemma is where the $\mu$-ergodicity is necessary for the proof of Proposition \ref{prop}.

\begin{lemma}\label{lerg}
If \  $x \in  \mho( r_i, \rho_l)$, $j \geq j(i, l)$, $k= k(i, l)$ and $\mu(P^{(k)}_j(x)) \geq \delta \mu(\hat P^{(k)}_j(x))$, we have:
$$
\mu(\hat P^{(k)}_j(x) \cap \mho) \geq \delta \mu(\hat P^{(k)}_j(x)).
$$
\end{lemma}
\begin{proof}
As $\mu$ is ergodic, for $\mu-$typical $y \in M$, we have that
$$
\mu(\hat P^{(k)}_j(x) \cap \mho)= \lim_{n \to +\infty} \frac 1 n \#\{1 \leq t \leq n; f^t(y) \in \hat P^{(k)}_j(x) \cap \mho \},
$$
and
$$
\mu(P^{(k)}_j(x))= \lim_{n \to +\infty} \frac 1 n \#\{1 \leq t \leq n; f^t(y) \in P^{(k)}_j(x) \}.
$$
By the definition of $\mho(r_i, \rho_l)$, between any pair of natural numbers $n_1$ and $n_2$ such that $f^{n_1}(y), f^{n_2}(y) \in P^{(k)}_j(x) \subset B'_r(x)$, there exists $n_1\leq t_1< n_2$, such that $f^{t_1}(y) \in (\ov{B'_{\rho_l r}(x)} \cap \mho) \subset (\hat P^{(k)}_j(x) \cap \mho)$. This implies that
$$
\#\{1 \leq t \leq n; f^t(y) \in \hat P^{(k)}_j(x) \cap \mho \} \geq
$$
$$
\#\{1 \leq t \leq n; f^t(y) \in P^{(k)}_j(x) \}- 1.
$$
Hence
$$
\mu(\hat P^{(k)}_j(x) \cap \mho) \geq \mu(P^{(k)}_j(x) ) \geq \delta \mu(\hat P^{(k)}_j(x) ).
$$

\end{proof}

Now define $ \La_\delta^0(i, l) $, for $\delta> 0$,  as the set of
points $x \in T^s$ such that for $k= k(i, l)$, we have
$$
\mu(P_j^{(k)}(x)) \geq \delta \mu(\hat P_j^{(k)}(x))
 \quad\text{ and }\quad
\mu(P_j^{(k)}(x)) \geq  \delta \mu(\tilde P_j^{(k)}(x)),
$$
 for an infinite
sequence $\vsi(x)$ of values of  $j$, and set
 $\La_\delta(i, l):= \La_\delta^0(i, l) \cap \mho(r_i, \rho_l)$.
The next lemma, a kind of Vitali's covering lemma, will be useful to estimate
the measure of $\mho^c$.

\begin{lemma} \label{levit}
Given a neighborhood $V$ of \  $\mho^c \cap
\La_\delta(i, l)$, there exist sequences $x_q \in \mho^c \cap \La_\delta(i, l)$, $(j_q), j_q \in \vsi(x_q)
\subset \natural$, $q= 1, 2, \dots$, such that
\begin{enumerate}
\item The sets $\hat P^{(k)}_{j_q}(x_q)$, $q \in \natural$ are disjoint and contained in $V$;

\item
$
\mu\big((\mho^c \cap \La_\delta(i, l)) \setminus \cup_{q \in \natural} \hat P^{(k)}_{j_q}(x_q)\big)= 0.
$
\end{enumerate}
\end{lemma}
\begin{proof}

By standard measure theoretical arguments, a translation $\tau: T^s
\to T^s$ can be found in such way that
$$
\mu(\tau(\cup \{\partial \hat A; A \in \SP^{(k)}_j, k \geq 1, j \geq
1\}))= 0;
$$
where $\partial \hat A$ is the boundary of $\hat A \in \hat
\SP^{(k)}_j$.
Denoting by $\SF$ the family of sets $P^{(k)}_j(x)$ with $x \in
\Lambda_\delta(i, l) \cap \mho^c$ and $j \in
\vsi(x)$. Take a sequence $A_u \in \SF$ satisfying:
\begin{enumerate}
\item $\hat A_u \subset V$, $\forall u \in \natural$, and $\mu(\hat
A_u \cap \hat A_e)= 0$, $\forall 1 \leq e< u$.

\item
\label{itemc}
$ \diam(A_u) = \max\{\diam(A); \hat A \subset V \text{ and }\mu(\hat
A \cap \hat A_e)= 0, \forall 1 \leq e< u\} $.
\end{enumerate}
Such properties imply that $\lim_{u \to + \infty} \diam(A_u)= 0$ and
\begin{equation}
\sum_{u} \mu(A_u)= \mu(\cup_u A_u) \leq 1. \label{eqsum}
\end{equation}

We claim that for $N \geq 1$
\begin{equation}
\big( \Lambda_\delta(i, l) \cap \mho^c \big) \setminus \cup_{u=
1}^N \ov {\hat A_u} \subset \cup_{u> N} \tilde A_u. \label{eqar}
\end{equation}

In fact, if $x \in \big( \Lambda_\delta(i, l) \cap \mho^c \big) \setminus \cup_{u= 1}^N \ov {\hat A_u}$, there exist $A \in \SF$ with $x \in A$ and
$
\hat A \cap (\cup_{u= 1}^N \ov{\hat A_u})= \emptyset.
$
Take $N_1> N$ such that $\hat A \cap \hat A_u= \emptyset$, $\forall 1 \leq u < N_1$
and $\hat A \cap \hat A_{N_1} \neq \emptyset$.
By item (\ref{itemc}) above, it follows that $\diam(\hat A) \leq \diam(\hat A_{N_1})$.
This implies that $\hat A \subset \tilde A_{N_1}$ and then
$$
x \in A \subset \tilde A_{N_1} \subset  \cup_{u> N} \tilde A_u,
$$
which concludes the proof of equation \ref{eqar}.
By such equation and our assumption that partition elements borders have zero measure,
we obtain that
$$
\mu\Big(
\big( \Lambda_\delta(i, l) \cap \mho^c \big) \setminus \cup_{u=
1}^N  {\hat A_u}
 \Big)=
\mu\Big(
\big( \Lambda_\delta(i, l) \cap \mho^c \big) \setminus \cup_{u=
1}^N \ov {\hat A_u}
 \Big) \leq
$$
$$
\mu\big(\cup_{u> N} \tilde A_u\big) \leq \sum_{u > N} \mu(\tilde A_u) \leq
\delta^{-1} \sum_{u > N} \mu(A_u).
$$
Due to eq. (\ref{eqsum}) the tail sum above goes to zero as $N \to +\infty$,
which implies
the lemma.

\end{proof}

Lemmas  \ref{lerg} and \ref{levit} are the key ingredients in the

\begin{theorem} \label{teou}
Let $M$ be a compact Riemannian manifold and let $f:M \to M$ be a measurable Borelian map
such that $\SM_1(f) \neq \emptyset$.
Let $\mho \subset M$ and $\mho( r, \rho)$ be Borelian subsets of $M$,
where $r> 0$, $\rho> 1$.

Suppose that the points $x \in \mho(r, \rho)$ have the following
property: if $y \in B_{r'}(x)$ for some $0< r' \leq r$ and $f^t(y)
\in B_{r'}(x)$, for some $t> 0$ then there exist $0 \leq t_1 \leq
t$, such that $f^{t_1}(y) \in \ov{B_{\rho r'}(x)} \cap \mho$.
Suppose also that { $r_i> 0$, $\rho_l > 1$ are two monotone
sequences} converging respectively to 0 and $+\infty$, such that
$\cup_{i, l} \mho(r_i, \rho_l)= M$. Then, $\mho$ has total
probability with respect to the map $f$.

\end{theorem}
\begin{proof}
Consider  $\La_\delta^0(i, l)$ and $\La_\delta(i, l)= \La_\delta^0(i, l) \cap \mho(r_i, \rho_l)$
the same sets defined above in our text.
By Lemma \ref{leprp}, this implies that
$$
\mu(\La_\delta^0(i, l)) \geq 1- \delta (k^s + 3^s k^s).
$$
Since the last inequality implies that
$
 \cup_{n= 1} ^{+\infty} \La_{1/n}(i, l)= \mho(r_i, \rho_l) (\!\mod 0)
$
it is enough to prove that
$$
\mu(\mho^c \cap \La_\delta(i, l))= 0, \forall 0< \delta< 1.
$$
In fact this implies that $\mho \supset \mho( r_i, \rho_l) \mod(0)$ and, consequently,
$$
 \mho \supset \big(\cup_{i, l} \mho( r_i, \rho_l)\big) = M \mod(0).
$$
We will then have
$\mu(\mho)= \mu(M)= 1$, and the proof of Theorem \ref{teou} will be completed.

Fix $(i, l)$ and $\delta> 0$. Let $V$ a neighborhood of  $\La_\delta(i, l) \cap \mho^c$.
 By lemmas
\ref{lerg} and \ref{levit} it follows that
$$
\mu(V) \geq \sum_q \mu\big(\hat  P^{(k)}_{j_q}(x_q)  \big) \geq
\frac{1}{1- \delta} \sum_q \mu\big(\hat  P^{(k)}_{j_q}(x_q) \cap \mho^c \big)=
$$
$$
\frac{1}{1- \delta}  \mu\Big(\big(\cup_q \hat  P^{(k)}_{j_q}(x_q) \big)\cap \mho^c \Big) \geq
\frac{1}{1- \delta}  \mu\big( \Lambda_\delta(i, l)\cap \mho^c \big).
$$
But if $\mu\big( \Lambda_\delta(i, l)\cap \mho^c \big)> 0$, one can take
$V$ satisfying
$$
\mu(V) < \frac{1}{1- \delta}  \mu\big( \Lambda_\delta(i, l)\cap \mho^c \big),
$$
contradicting the last inequality.
Hence $\mu\big( \Lambda_\delta(n, m)\cap \mho^c \big)= 0$.

\end{proof}

Now, let us finish the proof of Proposition \ref{prop} (which implies Theorem \ref{theop}).

\begin{proof}[Proof of Proposition \ref{prop}]


Define  $\Sigma(\SU, \epsilon, r, \rho)$,
where $r> 0$, $\rho> 1$,
as the set of
points $x \in M$ such that  if $y \in B_{r'}(x)$ for some $0< r'
\leq r$ and $f^t(y) \in B_{r'}(t)$, for some $t> 0$ then there exist $0 \leq t_1< t_2
\leq t$, $g \in \SU$ and $z \in M$ such that $g= f$ on $M \setminus
B_\epsilon(f, x)$,
$$
g^{t_2- t_1}(z)= z, d(g^j(z), f^j(f^{t_1}(y)) \leq \epsilon, \forall
0 \leq j \leq t_2- t_1,
$$
and
$$
 f^{t_1}(y) \in \ov{B_{\rho r'}(x)}.
$$
In particular, $f^{t_1}(y) \in \Sigma(\SU, \epsilon)$.
It is easy to see that  $\Sigma(\SU, \epsilon, r, \rho)$ is a Borelian set.
Again, let $r_i> 0$ and $\rho_l> 1$ to be two monotone sequences converging respectively to $0$
 and $+\infty$.
We note that
Lemma \ref{lesomb} implies that
$$
M= \cup_{i \geq 1} \cup_{l \geq 1} \Sigma(\SU, \epsilon, r_i, \rho_l),
$$
for every neighborhood $\SU$ of $f$ and every $\epsilon> 0$.
So, taking $\mho= \Sigma(\SU, \epsilon)$ and $\mho(r_i, \rho_l)= \Sigma(\SU, \epsilon, r_i, \rho_l)$ in Theorem \ref{teou}, we conclude that $\mu(\mho)= \mu(\Sigma(\SU, \epsilon))= 1$ for
all $f-$ergodic probability. By Ergodic Decomposition Theorem, this implies that $\Sigma(\SU, \epsilon)$ has total probability.

\end{proof}

So far, we have proven the raw version of Ergodic Closing Lemma for
Endomorphisms (Theorem \ref{theop}). The next lemma will be used in the
proof of the residual version of Ergodic Closing Lemma. We denote by
$\SM(M)$ the set of probabilities on $M$ endowed with the weak-*
topology.

\begin{lemma} \label{leres}
Let $f:M \to M$ be an endomorphism. Suppose that, for $x$ in a total
probability set $S \subset M$, given $\epsilon > 0$ and a
neighborhood  $\SU$ of $f$, there exists $g_{x,  \epsilon} \in \SU$
with a periodic point $p= p(x, \epsilon)$ which $\epsilon-$shadows
$x$. Then, given any ergodic measure $\mu \in \SM_1(f)$, there are
$g_k \to f$ and $g_k-$periodic points $p_k$ such that $\mu$
 is the limit of the sequence $(\mu_k)$ of $g_k-$ergodic measures respectivelly supported in
the orbit of $p_k$. Moreover, each $p_k$ can be taken to be a
hyperbolic periodic point for $g_k$.
\end{lemma}
\begin{proof}

Let us consider an $f$-ergodic probability $\mu$. We suppose,
without loss of generality, that $\mu$ is not supported in a
periodic orbit, otherwise there is nothing to prove. For a
$\mu-$typical point $x\in M$, we can assume that $x$ is recurrent
(by Poincar\'e's Recurrence Theorem), has the shadowing property as
in lemma's statement, and that
\begin{equation}
\frac 1 n \sum_{j= 0}^{n-1} \delta_{f^j(x)} \to_{_{\text{weak}-*}}
\mu, \label{eqfrac}
\end{equation}
as $n \to +\infty$. In the last claim we made use of the Ergodic Decomposition Theorem.

Set $\epsilon_1= 1$ and $n_k> 0$ as the first return time of the
orbit of $x$ to $B(x, \epsilon_k)$, where $\epsilon_{k+ 1}:=
d(f^{n_{k}}(x), x)/2$, $\forall k \geq 1$.

Therefore, $n_k \to +\infty$ as $k \to +\infty$. By hypothesis, one
can take  a sequence of $g_k:= g_{x, \epsilon_k}$, with $g_k \to f$,
exhibiting  $g_k-$periodic points $(p_k)$ such that each $p_k$
$\epsilon_k/3$-shadows the orbit of $x$. In particular, the period
$t_{k+1}$ of $p_{k+1}$ is, at least, $n_k$ (otherwise, the orbit of
$x$ would return to $B(x, \epsilon_k)$ before $n_k$). So, $t_{k+1}
\geq n_k$ implies that $t_k \to +\infty$ as $k \to +\infty$, and (up
to take a subsequence) we can suppose that $t_k$ are distinct. Note
that slightly perturbing $g_k$ in the neighborhood of $p_k$, we can
suppose that $p_k$ is hyperbolic. Set $\mu_k$ as the $g_k$-ergodic
probability supported in the orbit of $p_k$. We will show that
$\mu_k \to_{_{\text{weak}-*}} \mu$ as $k \to +\infty$. From equation
\ref{eqfrac} we have that
$$
\nu_k= \frac 1 {t_k} \sum_{j= 0}^{t_k-1} \delta_{f^j(x)}
\to_{_{\text{weak}-*}} \mu,
$$
as $k \to +\infty$. Let $\alpha> 0$ and $\{\vr_1, \dots,\vr_s\}
\subset C^0(M)$ be given. All we need to see is that there exists
$k_0 \in \natural$ such that $\mu_k$ belongs to the neighborhood
$$
V_{\vr_1, \dots, \vr_s; \alpha} := \{ \nu \in \SM(M); |\int \vr_i
d\nu - \int \vr_i d\mu|< \alpha, \forall i= 1, \dots, s\},
$$
forall $k \geq k_0$.
 In fact, as $\vr_i$, $i= 1, \dots, s$ are uniformly continuous,
 take $\epsilon> 0$ such that $|\vr_i(y)- \vr_i(z)|<
\alpha/2$, $\forall i= 1, \dots, s$, $\forall y, z \in M$ such
that $d(y, z)< \epsilon$. Then, take $k_0$ such that $\epsilon_k<
\epsilon/2$,and $|\int \vr_i d\nu_k - \int \vr_i d\mu | < \alpha/2$,
$\forall k \geq k_0$, $\forall i= 1, \dots, s$. We conclude that
$$
|\int \vr_i d\mu_k - \int \vr_i d\mu| \leq |\int \vr_i d\mu_k - \int
\vr_i d\nu_k| + |\int \vr_i d\nu_k - \int \vr_i d\mu| <
$$
$$
\frac 1 {t_k} \sum_{j= 0}^{t_k-1} |\vr_i(f^j(x))- \vr_i(g_k^j(p_k))|
+ \alpha/2 \leq \alpha, \forall i= 1, \dots, s;
$$
which implies the lemma.

\end{proof}

Now, we proceed with the proof of Theorem \ref{theo6}, by deriving it
from Theorem \ref{theop} and Lemma \ref{leres} above. The arguments here
are basically the same as in Theorem 4.2 in \cite{ABC}.

\vspace{.1cm}
\begin{proof}[ Proof of Theorem \ref{theo6}]

 For  $m \in \natural$ fixed, by standard transversality arguments,
the collection $\SK_m$ of  endomorphisms $f$ such that all periodic
points of $f$,  with period up to $m$ are hyperbolic is an open and
dense subset  $\NE(M)$.
So $\hat \SR:= \cap_{m= 1}^{+\infty} \SK_m$ is a residual set. Let
the set of probabilities $\SM(M)$ on $M$ to be endowed with the
weak-* topology and let $\kappa$ be the collection of compact
subsets of $\SM(M)$ endowed with Hausdorff distance. Given $f \in
\SR$, denote by $\SM_{per}(f)$ the set of $f-$ergodic measures
supported in $f-$periodic orbits.
Set $\Upsilon: \hat \SR \to \kappa$ given by
$$
\Upsilon(f)= \ov{\SM_{per}(f)}
$$
 Due to the
robustness of hyperbolic periodic points, such $\Upsilon$ is lower
semicontinuous. This implies that there is a residual subset $\SR
\subset \hat\SR$ whose elements are continuity points for
$\Upsilon$.

From now on, let $f \in \SR$. Let us prove that $\SM_1(f)$ is the
closed convex hull of $f-$ergodic measures supported in $f-$periodic
orbits. By Ergodic Decomposition Theorem, all we need to prove is
that any $f-$ergodic measure $\mu$ is in $\ov{\SM_{per}(f) }$.
By Lemma \ref{leres}, such measure $\mu$ is accumulated by $\mu_k
\in \SM_{per}(g_k)$, where $g_k \to f$ as $k \to +\infty$. As $\hat
\SR$ is residual, by means of a slight perturbation, we can suppose
that $g_k \in \hat \SR$ (as we construct $p_k$ to be hyperbolic in
the proof of that Lemma \ref{leres}, such $p_k$ persist under any
sufficiently small perturbation).
Since $f$ is a continuity point for $\Upsilon$, we have that
$\ov{\SM_{per}(g_k)} \to \ov{\SM_{per}(f)}$ as $k \to +\infty$, and
this implies that $\mu \in \ov{\SM_{per}(f)}$.

Therefore, $\ov{\SM_{per}(f)}$ contains all $f-$ergodic measures,
and by Ergodic Decomposition Theorem, we conclude that $\SM_1(f)$ is
the closed convex hull of $f$-ergodic measures supported in periodic
orbits.

\end{proof}

\end{document}